\documentclass[a4paper,12pt,pdf]{amsart}
\usepackage{enumerate, amsmath, amsfonts, amssymb, amsthm, thmtools, wasysym, graphics, graphicx, xcolor, frcursive,comment,bbm}

\usepackage{etex}

\usepackage[a4paper, margin=2.5cm]{geometry}

\usepackage{microtype}

\usepackage[hidelinks]{hyperref}

\usepackage{datetime}

\newdateformat{mydate}{\THEDAY~\monthname[\THEMONTH]~\THEYEAR}

\usepackage{calc}

\usepackage{tikz}
\usetikzlibrary{matrix}
\usetikzlibrary{shapes.geometric, positioning, calc, patterns}

\definecolor{darkblue}{rgb}{0.0,0,0.7} 
 
\definecolor{darkred}{rgb}{0.7,0,0} 
\usepackage{hyperref}
\usepackage[all]{xy}
\usepackage[T1]{fontenc}

\usepackage{caption,lipsum}
\usepackage{graphicx}                  
\usepackage{pstricks,pst-plot,pst-text,pst-tree,pst-eps,pst-fill,pst-node,pst-math}
\usepackage{setspace}
\usepackage{multicol}
\usepackage{tikz-cd}

\usetikzlibrary{positioning}

\newcommand{\darkred}{\color{darkred}} 
\newcommand{\defn}[1]{\emph{\darkred #1}}

\usetikzlibrary{arrows.meta}

\usepackage{etex}

\newtheorem{theorem}{Theorem}[section]
\newtheorem{prop}[theorem]{Proposition}
\newtheorem{lemma}[theorem]{Lemma}
\newtheorem{cor}[theorem]{Corollary}

\theoremstyle{definition}
\newtheorem{definition}[theorem]{Definition}
\newtheorem{rmq}[theorem]{Remark}

\numberwithin{equation}{section}

\title[Trickle groups on a complete graph are linear]{Trickle groups on a complete graph are linear}

\author{Thomas Gobet}

\address{Université Clermont Auvergne, LMBP, UMR 6620 (CNRS), Campus des Cézeaux, 3 place Vasarely, TSA 60026, CS 60026, 63178 Aubière cedex, France}

\begin{document}
	\maketitle
	
	\begin{abstract}
	Trickle groups were recently introduced by Bellingeri, Godelle and Paris~\cite{trickle}, as generalizations of both cactus groups and right-angled Coxeter and Artin groups. We construct a representation of any trickle group $\mathrm{Tr}(\Gamma)$ with complete graph $\Gamma$ on a space having a basis in bijection with the vertex set $V(\Gamma)$, and show that it is faithful. This in particular shows that trickle Garside groups are linear. It also allows one to realize every finite trickle group as a subgroup of $\mathrm{GL}_{|V(\Gamma)|}(\mathbb{C})$. We do not use the results from~\cite{trickle}; instead we rather show that several properties of trickle groups from~\cite{trickle} (and also some new ones) can easily be proven in the particular case of a complete graph using our representation. This includes the normal forms, the embedding of the trickle monoid inside the group in the case where $\mu(x)=\infty$ for all $x\in V(\Gamma)$, and the absence of torsion of the latter. In this last case we give a second proof of faithfulness in the spirit of Krammer-Hée, which explains how to read Garside-theoretic information from the representation. 
	\end{abstract}

	\section{Introduction}
	
	Trickle groups were introduced in 2024 by Bellingeri, Godelle and Paris~\cite{trickle}. A trickle group $\mathrm{Tr}(\Gamma)$ is defined by generators and relations, with relations of the form $xy=yz$ and $x^{\mu(x)} = 1$, where $x,y,z$ are generators; the presentation is defined using combinatorial data encapsulated by a graph $\Gamma$ called a \textit{trickle graph}. It is a simplicial graph $\Gamma$, endowed with a partial order $\leq$ on the set $V(\Gamma)$ of vertices of $\Gamma$, a label $\mu(x) \in \{2, 3, \dots\}\cup\{\infty\}$ for each $x\in V(\Gamma)$, and for each $x\in V(\Gamma)$, an automorphism $\varphi_x$ of the star graph of $x$ in $\Gamma$, satisfying certain compatibility conditions (see Section~\ref{sec:prelim} below or~\cite{trickle} for the precise definition). While these compatibility conditions are strong, there are many quite different families of groups fitting into this setting, and coming from different areas of mathematics: this includes cactus groups, right-angled Coxeter and Artin groups, ordered quandle groups, but also Thompson's group $F$, etc--we refer the reader to~\cite[Section 3]{trickle} for other examples. 
	
	Trickle groups have a solvable word problem. In~\cite{trickle}, an algorithm, called the \textit{trickle algorithm}, allows one to calculate normal forms for elements of such groups. It is based on a terminating and confluent rewriting system. A second algorithm in the spirit of Tits is given, generalizing the solution to the word problem in both right-angled Coxeter and Artin groups. While the relations $x^{\mu(x)} = 1$ show that these groups have torsion in general, in the case where $\mu(x) = \infty$ for all $x\in V(\Gamma)$, the defining presentation is positive, and one can consider the \textit{trickle monoid} $\mathrm{Tr}^+(\Gamma)$ defined by the same presentation. Bellingeri, Godelle and Paris show that $\mathrm{Tr}^+(\Gamma)$ is a Garside monoid if and only if $V(\Gamma)$ is finite and the graph $\Gamma$ is complete, in which case $\mathrm{Tr}(\Gamma)$ is thus a Garside group. Garside monoids are monoids satisfying certain divisibility conditions which ensure that they embed in their groups of fractions (which have the same presentation), and that the latter admit a solution to most combinatorial questions: word problem, determination of the center, absence of torsion, etc. -- we refer the reader to~\cite{DP, Garside} for more on Garside monoids and groups. 
	
	In this case, that is, if $\mu \equiv \infty$, $V(\Gamma)$ is finite, and $\Gamma$ is complete, this thus yields a new family of Garside groups. In addition, considering the same trickle graph but with finite values of $\mu$ yields a trickle group which is always finite. This is a very particular behaviour, showing that Garside trickle groups are in some sense combinatorially much simpler than other families of Garside groups: for Artin groups of spherical type and torus knot groups (more generally circular groups), two of the easiest families of Garside groups, there are only finitely many finite quotients that one can obtain from a group in these families by adding torsion on the generators~\cite{Gobet_toric}. 
	
	Another property that many families of Garside groups share is that they are linear, hence residually finite. While it is not true in general that every Garside group is linear~\cite{HH}, it seems natural to investigate whether given families of Garside groups are linear. One aim of this note is to show that Garside trickle groups are linear, and the representation used to show this property is built in a particularly simple way from the automorphisms $\varphi_x$. We show more generally that, for every trickle graph $\Gamma$ with a complete graph (which can have infinitely many vertices, and finite values of $\mu$), there is a faithful representation of $\mathrm{Tr}(\Gamma)$ on a space having a basis indexed by the vertex set $V(\Gamma)$; it is thus infinite-dimensional when $V(\Gamma)$ is infinite. This construction is somewhat reminiscent of the Tits' representation of Coxeter groups. Since a trickle group is finite if and only if $\Gamma$ is complete, $V(\Gamma)$ is finite, and $\mu(x) < \infty$ for all $x\in V(\Gamma)$, this also allows one to realize every finite trickle group as a subgroup of $\mathrm{GL}_{|V(\Gamma)|}(\mathbb{C})$. 
	
	We also show that our representation can be used to easily recover several results from~\cite{trickle} for trickle groups on a complete graph. This include the normal forms, hence the solution to the word problem, and in the case where $\mu \equiv \infty$, the embedding of the trickle monoid inside the group. While this setting is much more restrictive than the general case, the corresponding results in~\cite{trickle} are obtained as a corollary of the trickle algorithm, which is very involved. Here we only use elementary notions, in particular we do not use the results from~\cite{trickle}, and this note is mostly self-contained. 
	
	Our approach shows that the normal forms of elements of the concerned trickle groups can easily be recovered from their matrix. In short, one shows that it is possible to write elements under a certain form, and then that the matrix reads this form, showing that it is a normal form. Since our setting also covers the case where the trickle groups are Garside, in Section~\ref{sec_alt} we give a second proof of faithfulness for the trickle monoid on a complete graph (implying faithfulness for the group in the case where the graph is finite), which is in the spirit of Krammer's proof of the linearity of braid groups, or rather Hée's simplified approach~\cite{Hee}. We show that we can easily recover atoms left- or right-dividing an element from the matrix of the element. This yields a new family of Garside groups for which Krammer's approach can be successfully applied to show linearity. This approach was already used independently by Cohen-Wales~\cite{CW} and Digne~\cite{Dig} for Artin groups of spherical type, and by the author for torus knot groups~\cite{Gobet_Burau}. It was also used by Paris to show that Artin monoids (which are not Garside in general) embed in their groups~\cite{Paris_inject}. 
	
	Our representation does not admit an obvious generalization to all trickle groups. It would be interesting to be able to extend it to all trickle graphs, or to right-angled Artin groups, and try to show that it is faithful. 
	
 \medskip
	
	\textbf{Funding statement}. The author is partially supported by the ANR (Agence Nationale de la Recherche) project CaGeT (ANR-25-CE40-4162).
	
	\textbf{AI statement}. The author did not use AI to assist with any kind of aspect of this manuscript.

	\section{Trickle groups}\label{sec:prelim}
	
	We recall the definition of trickle graphs and groups from~\cite{trickle}. If $\Gamma$ is a simplicial graph, we denote by $V(\Gamma)$ its set of vertices and by $E(\Gamma)$ its set of edges. Given a vertex $x\in V(\Gamma)$, the \defn{star} of $x$ is the full subgraph $\mathrm{star}_x(\Gamma)$ of $\Gamma$ spanned by $\{ y\in V(\Gamma) \mid \{x,y\}\in E(\Gamma)\} \cup \{x\}$.  
	
	\begin{definition}[Trickle graph]\label{def_graph}
	A \defn{trickle graph} is a quadruple $(\Gamma, \leq, \mu, (\varphi_x)_{x\in V(\Gamma)})$ where $\Gamma$ is a simplicial graph, $\leq$ is a partial order on $V(\Gamma)$, $\mu$ is a function $V(\Gamma) \longrightarrow \mathbb{Z}_{\geq 2} \cup\{ \infty\}$, and for all $x\in V(\Gamma)$, $\varphi_x$ is an automorphism of graph of $\mathrm{star}_x(\Gamma)$, such that 
	\begin{itemize}
	\item[(a)] For all $x, y\in V(\Gamma)$, if $x < y$, then $\{x,y\} \in E(\Gamma)$, 
	\item[(b)] For all $x, y, z\in V(\Gamma)$, if $\{x,y\}\in E(\Gamma)$ but $x$ and $y$ are not comparable for $\leq$, and $z\leq x$, then $z$ and $y$ are not comparable for $\leq$,  
	\item[(c)] For all $x\in V(\Gamma)$ and $y, z\in \mathrm{star}_x(\Gamma)$, we have $y \leq z$ if and only if $\varphi_x(y) \leq \varphi_x(z)$, 
	\item[(d)] For all $x\in V(\Gamma)$, if $y \in \mathrm{star}_x(\Gamma)$ and $\varphi_x(y) \neq y$, then $y<x$, 
	\item[(e)] For all $x\in V(\Gamma)$, if $\mu(x)$ is finite, then $\varphi_x$ has finite order, and its order divides $\mu(x)$, 
	\item[(f)] For all $x\in V(\Gamma)$ and $y\in \mathrm{star}_x(\Gamma)$, we have $\mu(x) = \mu(\varphi_x(y))$, 
	\item[(g)] For all $x, y, z\in V(\Gamma)$ such that $z < y < x$, setting $y' = \varphi_x(y)$, we have $$(\varphi_x \circ \varphi_y) (z) = (\varphi_{y'} \circ \varphi_x) (z).$$ 
	\end{itemize}
	\end{definition}
	
	If $x, y\in V(\Gamma)$ are not comparable for $\leq$, we shall simply write $x \Vert y$. 
	
	\begin{definition}[Trickle group]
	Let $\Gamma = (\Gamma, \leq, \mu, (\varphi_x)_{x\in V(\Gamma)})$ be a trickle graph. The \defn{trickle group} attached to the trickle graph $\Gamma$ is the group $\mathrm{Tr}(\Gamma)$ defined by the presentation $$\langle x\in V(\Gamma) \mid x^{\mu(x)} = 1, \ \forall x\in V(\Gamma)\text{~with~} \mu(x)< \infty, \ \varphi_x(y) x = \varphi_y(x) y, \ \forall \{x, y\}\in E(\Gamma)\rangle.$$
	\end{definition}
	
Examples of trickle groups include right-angled Coxeter and Artin groups, which correspond to the case where $\leq$ is the trivial order, i.e., any two distinct elements from $V(\Gamma)$ are never comparable, and $\mu(x) = 2$ (respectively $\mu(x)= \infty$) for all $x\in V(\Gamma)$. Other examples include cactus groups, Thompson's group $F$, etc.-- we refer the reader to~\cite{trickle} for more examples of trickle groups.

In the next sections, we will restrict ourselves to those trickle groups attached to complete graphs. In this case, we have $\mathrm{star}_x(\Gamma) = \Gamma$ for all $x\in V(\Gamma)$, hence $\varphi_x$ induces an automorphism of poset of $(V(\Gamma), \leq)$ for all $x\in V(\Gamma)$. The key point will be that, in this case, $x\mapsto \varphi_x$ defines a representation of $\mathrm{Tr}(\Gamma)$. We shall linearize and twist it in Section~\ref{sec_rep} below to get a faithful representation of $\mathrm{Tr}(\Gamma)$.

	\section{A faithful representation of trickle groups with a complete graph}\label{sec_rep}

	Let $\Gamma=(\Gamma, \leq, \mu, (\varphi_x)_{x\in V(\Gamma)})$ be a trickle graph. We assume that $\Gamma$ is a complete graph. In this section, we shall construct a faithful representation of the trickle group $\mathrm{Tr}(\Gamma)$. 
	
	\subsection{Construction of the representation}
	
	Since $\Gamma$ is complete, for all $x\in V(\Gamma)$ we have $\mathrm{star}_x(\Gamma) = V(\Gamma)$. We thus have that $\varphi_x$ is an automorphism of poset of $(V(\Gamma), \leq)$ for all $x\in V(\Gamma)$. Let $G_\Gamma$ denote the subgroup of the group of automorphisms of posets of $(V(\Gamma), \leq)$ generated by $(\varphi_x)_{x\in V(\Gamma)}$. Let $O_\Gamma$ denote the set of orbits of $G_\Gamma$ on $V(\Gamma)$. Given $x\in V(\Gamma)$ we denote by $[x]$ its $G_\Gamma$-orbit. Note that by condition~$(f)$ in Definition~\ref{def_graph}, the labeling map $\mu$ is constant on $G_\Gamma$-orbits. We denote by $O_\Gamma^\infty$ the set of orbits $[x]$ such that $\mu(x) = \infty$.  
		
	For every orbit $[x]$ in $O_\Gamma$ such that $\mu(x) = \infty$, we let $q_{[x]}$ be an indeterminate. Let $R$ denote the ring of Laurent polynomials in the indeterminates $q_{[x]}$, $[x] \in O_\Gamma^\infty$ with coefficients in $\mathbb{C}$. For $x\in V(\Gamma)$ we set $n_x = q_{[x]}$ if $\mu(x) = \infty$ and $n_x = e^{2 i \pi / \mu(x)}$ otherwise. 
	
Let $V$ the the free $R$-module with a basis $\mathcal{B}=\{e_x \mid x\in V(\Gamma)\}$. We will turn $V$ into a representation of $\mathrm{Tr}(\Gamma)$. To this end we define the element $\rho_x$ of $\mathrm{GL}(V)$ attached to a generator $x\in V(\Gamma)$ by its action on the elements of $\mathcal{B}$ as follows:
\[ \rho_x( e_y) = \begin{cases} n_x e_x \text{~if~} y=x, \\
e_{\varphi_x(y)}\text{~otherwise.}\end{cases}\]
In particular, if $y \not< x$, then $\varphi_x(y) = y$, hence $\rho_x(e_y) = e_y$. 

\begin{prop}\label{prop_rho}
The assignment $x \mapsto \rho_x$ extends to a group homomorphism $\rho : \mathrm{Tr}(\Gamma) \longrightarrow \mathrm{GL}(V)$, that is, a representation of $\mathrm{Tr}(\Gamma)$.  
\end{prop}

\begin{proof}
It suffices to show that the elements $(\rho_x)_{x\in V(\Gamma)}$ satisfy the defining relations of $\mathrm{Tr}(\Gamma)$. 

Let $x\in V(\Gamma)$ such that $n:=\mu(x) < \infty$. We show that $\rho_x^n = \mathrm{id}_V$. Since $\varphi_x$ is a bijection of $V(\Gamma)$ and $\varphi_x(x)=x$, for all $y \neq x$ and all $m\in \mathbb{Z}$ we have $\varphi_x^m(y) \neq x$. We thus have $\rho_x^n(e_x) = n_x^n e_x = e_x$, and for $y\neq x$, we have $\rho_x^n(e_y) = e_{\varphi_x^n(y)}$. Since the order of $\varphi_x$ divides $\mu(x) = n$, we thus have $\rho_x^n(e_y) = e_y$. This shows that $\rho_x^n = \mathrm{id}_V$.

Now let $x, y\in V(\Gamma)$ such that $\{x,y\}\in E(\Gamma)$. We need to show that $\rho_{\varphi_x(y)}\circ \rho_x = \rho_{\varphi_y(x)} \circ\rho_y$. We divide the proof into two cases, depending on whether $x < y$ or $x \Vert y$ (the case $y < x$ is obtained from the first one by reversing the roles of $x$ and $y$). 

\medskip

\textit{Case $x < y$}. Note that $\varphi_x(y)=y$. Let $z\in V(\Gamma)$. We need to show that $\rho_{y}\circ \rho_x(e_z)=\rho_{\varphi_y(x)} \circ\rho_y(e_z)$. Assume that $z = x$. We have $$\rho_y \circ \rho_x( e_z) = \rho_y ( n_x e_x) = n_x e_{\varphi_y(x)} = n_{\varphi_y(x)} e_{\varphi_y(x)} = \rho_{\varphi_y(x)}(e_{\varphi_y(x)}) =\rho_{\varphi_y(x)} \circ \rho_y (e_z).$$ Assume that $z = y$. First note that since $x < y$ and $\varphi_y(y)=y$, we have $\varphi_y(x) < y$, and hence $\varphi_{\varphi_y(x)}(y) = y$. We have   
$$\rho_y \circ \rho_x( e_z) = \rho_y ( e_y) = n_y e_{y} = n_{y} \rho_{\varphi_y(x)}( e_y) = \rho_{\varphi_y(x)}(n_y e_y) =\rho_{\varphi_y(x)} \circ \rho_y (e_z).$$

We can thus assume that $z\neq x$ and $z\neq y$. We treat the remaining cases, depending on whether $y < z$, $y \Vert z$, or $z < y$. 

If $y < z$, then $\varphi_x(z) = \varphi_y(z) = z$, and since $x < z$ we have $\varphi_y(x) < \varphi_y(z)=z$, yielding $\varphi_{\varphi_y(x)} (z) =z$. We thus have $\rho_{y}\circ \rho_x(e_z)=e_z=\rho_{\varphi_y(x)} \circ\rho_y(e_z)$.

If $y \Vert z$, then since $x < y$ and $\varphi_y(x) < y$, by condition~$(b)$ in Definition~\ref{def_graph} we have $x \Vert z$ and $\varphi_y(x) \Vert z$. This yields $\varphi_x(z) = \varphi_y(z) = \varphi_{\varphi_y(x)}(z) = z$, and hence $\rho_{y}\circ \rho_x(e_z)=e_z=\rho_{\varphi_y(x)} \circ\rho_y(e_z)$. 

If $x < z < y$, then $\varphi_y(x) < \varphi_y(z)$ hence $\varphi_{\varphi_y(x)}({\varphi_y(z)}) = \varphi_y(z)$, and we get $$\rho_y \circ \rho_x(e_z) = \rho_y(e_z) = e_{\varphi_y(z)} = \rho_{\varphi_y(x)}( e_{\varphi_y(z)}) = \rho_{\varphi_y(x)}\circ \rho_y (e_z).$$ 

If $z < y$ and $z \Vert x$, then $\varphi_y(z) \Vert \varphi_y(x)$, hence $\varphi_x(z)=z$ and $\varphi_{\varphi_y(x)}(\varphi_y(z)) = \varphi_y(z)$, and the calculation is the same as in the previous case. 

Finally, if $z < x < y$, then by condition~$(g)$ in Definition~\ref{def_graph} we have $\varphi_y \circ \varphi_x(z) = \varphi_{\varphi_y(x)} \circ \varphi_y (z)$. We thus have $$\rho_y \circ \rho_x(e_z) = \rho_y (e_{\varphi_x(z)}) = e_{\varphi_y \circ \varphi_x(z)} = e_{\varphi_{\varphi_y(x)} \circ \varphi_y (z)}=\rho_{\varphi_y(x)}(e _{\varphi_y}(z)) = \rho_{\varphi_y(x)} \circ \rho_y (e_z).$$

\textit{Case $x\Vert y$}. In this case we have $\varphi_y(x)=x$ and $\varphi_x(y)=y$. The relation to check is hence the commutation relation $xy=yx$. Assume that $z=x$. We then have $$\rho_y \circ \rho_x( e_z) = \rho_y (n_x e_x) = n_x e_x = \rho_x( e_x) = \rho_x \circ \rho_y (e_z),$$ and the case where $z=y$ is obtained by exchanging the roles of $x$ and $y$ which are symmetric here. 

We can thus assume that $z\neq x$ and $z\neq y$ and again, we treat the remaining cases depending on the position of $z$ with respect to $y$. 

If $ y < z$, then since $x \Vert y$ we either have $x < z$ or $ x \Vert z$. In both cases we have $\varphi_x(z) = z$. We get $$\rho_y \circ \rho_x (e_z) = \rho_y (e_z) = e_z = \rho_x( e_z) = \rho_x \circ \rho_y(e_z).$$

If $z < y$, then by condition~$(b)$ in Definition~\ref{def_graph} we have $z \Vert x$, hence $\varphi_y(z) \Vert \varphi_y(x)=x$. This yields $\varphi_x(\varphi_y(z))=\varphi_y(z)$. We get $$\rho_x \circ \rho_y (e_z) = \rho_x( e_{\varphi_y(z)}) = e_{\varphi_y(z)}=\rho_y(e_z)=\rho_y \circ \rho_x (e_z).$$

If $z \Vert y$, then we can assume that $z \Vert x$, otherwise we argue as in the first case if $z >x$ or as in the previous case if $z < x$, with the roles of $x$ and $y$ reversed. We then get $$\rho_y \circ \rho_x( e_z) = e_z = \rho_x \circ \rho_y (e_z).$$  
\end{proof}
	
	\subsection{Faithfulness}
	
	The aim of this subsection is to show that the representation built in Section~\ref{sec_rep} is faithful. We keep the assumptions made on $\Gamma$ in Section~\ref{sec_rep}. Given $g\in \mathrm{Tr}(\Gamma)$ we denote by $\rho_g$ its image under the homomorphism $\rho$ from Proposition~\ref{prop_rho}. 
	
	We fix a linear extension $\preceq$ of $\leq$, that is, a total order $\preceq$ on $V(\Gamma)$ such that for all $x, y\in V(\Gamma)$, $x \leq y \Rightarrow x \preceq y$. Given $g\in \mathrm{Tr}(\Gamma)$, we first establish that there exists a word representing $g$ of the form $x_1^{a_1} x_2^{a_2} \cdots x_k^{a_k}$, where $x_i \in V(\Gamma)$, $x_1 \succ x_2 \succ \cdots \succ x_k$, and $a_i \in \mathbb{Z} \backslash \{0\}$ (the identity corresponds to the case $k=0$). Note that this result is already known: it is a consequence of the trickle algorithm of~\cite[Section 4]{trickle}, which holds for arbitrary trickle groups, and yields a normal form of the above form for every element of $\mathrm{Tr}(\Gamma)$ in the case where the graph is complete. But since the trickle algorithm is quite involved, we will give an elementary proof of the existence and uniqueness of such a form, reproving Bellingeri-Godelle-Paris' result in this specific case using our representation. 
	
	\begin{lemma}\label{lemma_form}
	Let $g\in \mathrm{Tr}(\Gamma)$. There exists a word representing $g$ of the form $x_1^{a_1} x_2^{a_2} \cdots x_k^{a_k}$, where $x_i \in V(\Gamma)$, $x_1 \succ x_2 \succ \cdots \succ x_k$, and $a_i \in \mathbb{Z}$.
	\end{lemma}
	
	\begin{proof}
	We start from any word $y_1^{\varepsilon_1} y_2^{\varepsilon_2} \cdots y_n^{\varepsilon_n}$, where $y_i \in V(\Gamma)$ and $\varepsilon_i \in\{-1, 1\}$. We show by induction on $n$ that using the defining relations, we can transform this word into a word of the form ${x_1}^{\delta_1} x_2^{\delta_2} \cdots x_n^{\delta_n}$, where $\delta_i\in\{-1, 1\}$, and $x_1$ is the element of $\{y_1, y_2, \dots, y_n\}$ which is maximal for the partial order $\preceq$, and $x_1 \succeq x_2 \succeq \dots \succeq x_n$. If $n=0$ or $n=1$ the result is clear. Now assume that $n>1$ and that the result holds for words of length smaller than $n$. We show that we can move the first instance of $x_1$ in the word $y_1^{\varepsilon_1} y_2^{\varepsilon_2} \cdots y_n^{\varepsilon_n}$ to the left of the word using the defining relations, without changing the represented element $g$ nor the length of the word, and such that the letters other than $x_1$ in the resulting word are all smaller than or equal to $x_1$ for $\preceq$. If $y_1=x_1$ then we are done. Hence assume that $y_1\neq x_1$ and let $i> 1$ be the first integer such that $y_i =x_1$. By definition of $i$, since $x_1$ is maximal for $\prec$, we have $y_{i-1} \prec x_1$. We show that we can replace the portion $y_{i-1}^{\varepsilon_{i-1}} y_i^{\varepsilon_i}$ of our word by a word of the form $x_1^{\delta} y^{\delta'}$, with $y\in V(\Gamma)$, $y \prec x_1$, and $\delta, \delta'\in\{-1,1\}$. If $x_1 \Vert y_{i-1}$ then we have $x_1 y_{i-1} = y_{i-1} x_1$, and the result is obvious. We can thus assume that $y_{i-1} < x_1$. Assume that $\varepsilon_{i} = 1$. Let $z\in V(\Gamma)$ such that $\varphi_{x_1}(z) = y_{i-1}$. We have a defining relation of $\mathrm{Tr}(\Gamma)$ given by $\varphi_{x_1}(z) x_1 = \varphi_z(x_1) z$. Now since $y_{i-1} < x_1$, we have $\varphi_{x_1}^{-1}(y_{i-1}) < \varphi_{x_1}^{-1}(x_1) = x_1$, hence $z < x_1$, yielding $\varphi_z(x_1) = x_1$. The relation is thus $y_{i-1} x_1 = x_1 z$, and we can replace $y_{i-1} x_1$ by $x_1z$ (if $\varepsilon_{i-1}=1$) or $y_{i-1}^{-1} x_1$ by $x_1 z^{-1}$ (if $\varepsilon_{i-1}=-1$). In this way $x_1$ has been moved one position to the left, and the element $z$ is still smaller than $x_1$ for $\prec$. Now assume that $\varepsilon_i = -1$. We consider the relation $\varphi_{y_{i-1}}(x_1) y_{i-1} = \varphi_{x_1}(y_{i-1}) x_1$; since $y_{i-1} < x_1$ we have $\varphi_{y_{i-1}}(x_1)=x_1$ and $z:=\varphi_{x_1}(y_{i-1}) < x_1$. The relation is thus $x_1 y_{i-1} = z x_1$, and we can replace $y_{i-1} {x_1}^{-1}$ by ${x_1}^{-1} z$ (if $\varepsilon_{i-1}=1$) or $y_{i-1}^{-1} {x_1}^{-1}$ by ${x_1}^{-1} z^{-1}$ (if $\varepsilon_{i-1}=-1$). Again, we moved the instance of $x_1$ one position to the left with the required properties. Repeating, we move $x_1$ to the left of the word and get a word $x_1^{\delta_1} z_2^{\gamma_2} \cdots z_n^{\gamma_n}$, where $z_i \preceq x_1$ for all $i$ and $\delta_1, \gamma_i \in \{-1,1\}$. Applying the induction hypothesis to the word $z_2^{\gamma_2} \cdots z_n^{\gamma_n}$ yields the result.
	
	To conclude the proof we apply the procedure above to any word in $V(\Gamma) \cup V(\Gamma)^{-1}$ representing the element $g$. Note that this yields a word where the powers of generators are $\pm 1$, but $x_i$ and $x_{i+1}$ may be equal. Grouping the same generators yields the claimed form.   
  
	\end{proof}
	
	\begin{rmq}\label{rmq_monoid}
	In the case where $\mu(x) = \infty$ for all $x\in V(\Gamma)$, one defines the \defn{trickle monoid} $\mathrm{Tr}^+(\Gamma)$ as the monoid with the same presentation as $\mathrm{Tr}(\Gamma)$. The above Lemma also works in this case, i.e., for elements of $\mathrm{Tr}^+(\Gamma)$, with the condition $a_i \in \mathbb{Z}$ replaced by $a_i \geq 1$ for all $i$. Indeed, the same proof can be given, only using the case $\varepsilon_{i-1} = 1 = \varepsilon_i$. 
	\end{rmq}
	
	As a consequence, every element $g\in \mathrm{Tr}(\Gamma)$ can be represented by a word of the form $x_1^{a_1} x_2^{a_2} \cdots x_k^{a_k}$, where $x_i \in V(\Gamma)$, $x_1 \succ x_2 \succ \cdots \succ x_k$, and $ 1 \leq a_i \leq \mu(x_i)-1$ if $\mu(x_i) < \infty$, and $a_i \in \mathbb{Z} \backslash \{0\}$ if $\mu(x_i)= \infty$. Indeed, starting from the word obtained from Lemma~\ref{lemma_form}, we remove an $x_i^{a_i}$ if $a_i=0$ or if $\mu(x_i)< \infty$ and $a_i$ is divisible by $\mu(x_i)$, and otherwise (still if $\mu(x_i) < \infty$) we replace it by its unique representative in $\{1, 2, \dots, \mu(x_i)-1\}$. We keep its power if $\mu(x_i) = \infty$.

\begin{lemma}\label{lem_main}
Let $g\in \mathrm{Tr}(\Gamma)$ be represented in the form $x_1^{a_1} x_2^{a_2} \cdots x_k^{a_k}$, where $x_i \in V(\Gamma)$, $x_1 \succ x_2 \succ \cdots \succ x_k$, and $ 1 \leq a_i \leq \mu(x_i)-1$ if $\mu(x_i) < \infty$, and $a_i \in \mathbb{Z} \backslash \{0\}$ if $\mu(x_i)= \infty$. Let $\varphi = \varphi_{x_1}^{a_1} \circ \varphi_{x_2}^{a_2} \circ \cdots \circ \varphi_{x_k}^{a_k}$. We have 
$$ \rho_g(e_y) = \begin{cases} n_{x_i}^{a_i} e_{\varphi(y)} &\text{~if~} y = x_i, \\
	e_{\varphi(y)}&\text{~otherwise}.
	\end{cases}$$
\end{lemma}

\begin{proof}
We argue by induction on $k$. If $k=0$ then $g=1$ and $\varphi = \mathrm{id}$, and the claim is obvious. Assume that $k\geq 1$ and that the result hold by induction for $\ell < k$. Since $x_1 \succ x_i$ for all $i=2, \dots, k$, we have $\varphi_{x_i}(x_1) = x_1$ for all $i=1, \dots, k$, hence $\rho_g(e_{x_1}) = \rho_{x_1}^{a_1} (e_{x_1}) = n_{x_1}^{a_1} e_{x_1}$. Since $\varphi(x_1) = x_1$ we get the result in this case. Let $i\in \{2, \dots, k\}$. Setting $g' = x_2^{a_2} \cdots x_k^{a_k}$, by induction we have $\rho_{g'}(e_{x_i}) = n_{x_i}^{a_i} e_{\varphi'(x_i)}$, with $\varphi'=\varphi_{x_2}^{a_2} \circ \cdots \circ \varphi_{x_k}^{a_k}$. Now since $x_1 \not\leq x_2, \dots, x_k$, we have $\varphi'(x_1) = x_1$. Since $x_i \neq x_1$ and $\varphi'$ is a bijection, we have $\varphi'(x_i) \neq x_1$, and hence, $(\varphi_{x_i}^m \circ \varphi')(x_i) \neq x_1$ for all $m \geq 0$. This yields $$\rho_g(e_{x_i}) = \rho_{x_1}^{a_1} ( n_{x_i}^{a_i} e_{\varphi'(x_i)})= n_{x_i}^{a_i} e_{\varphi_{x_1}^{a_1} ( \varphi'(x_i))} = n_{x_i}^{a_i} e_{\varphi(x_i)}.$$
Now assume that $y\neq x_i$ for all $i=1, \dots, k$. By induction we have $\rho_{g'}(e_y) = e_{\varphi'(y)}$. Since $\varphi'(x_1) = x_1$ and $\varphi'$ is a bijection, we have $\varphi'(y) \neq x_1$. Arguing as in the previous case yields $\rho_g(e_y) = e_{\varphi(y)}$. 
\end{proof}

We now easily deduce our main result from the above lemma. The points $(2)$ to $(5)$ below in particular give new proofs of results from~\cite{trickle} (see Corollary 2.7 and its proof, Theorems 2.16 and 2.17) in the particular case where $\Gamma$ is complete.  

\begin{theorem}\label{main}
Let $\Gamma$ be trickle graph. Assume that $\Gamma$ is complete. 
\begin{enumerate}
	\item The representation $\rho$ is faithful, 
	\item Every element $g\in \mathrm{Tr}(\Gamma)$ can be written uniquely in the form $x_1^{a_1} x_2^{a_2} \cdots x_k^{a_k}$, where $x_i \in V(\Gamma)$, $x_1 \succ x_2 \succ \cdots \succ x_k$, and $ 1 \leq a_i \leq \mu(x_i)-1$ if $\mu(x_i) < \infty$, and $a_i \in \mathbb{Z} \backslash \{0\}$ if $\mu(x_i)= \infty$.
	\item The group $\mathrm{Tr}(\Gamma)$ is finite if and only if $\mu(x) < \infty$ for all $x\in V(\Gamma)$, in which case its order is $\prod_{x\in V(\Gamma)} \mu(x)$,
	\item In the case where $\mu(x)= \infty$ for all $x\in V(\Gamma)$, the monoid $\mathrm{Tr}(\Gamma)^+$ embeds into $\mathrm{Tr}(\Gamma)$. In particular, it is cancellative.
	\item In the case where $\mu(x)= \infty$ for all $x\in V(\Gamma)$, the group $\mathrm{Tr}(\Gamma)$ is torsion-free.  
\end{enumerate}
\end{theorem}

\begin{proof}
\begin{enumerate}
	\item Let $g\in \mathrm{Tr}(\Gamma)$, $g\neq 1$. We need to show that $\rho_g \neq \mathrm{id}_V$. Write $g$ in the form $x_1^{a_1} x_2^{a_2} \cdots x_k^{a_k}$, where $x_i \in V(\Gamma)$, $x_1 \succ x_2 \succ \cdots \succ x_k$, and $ 1 \leq a_i \leq \mu(x_i)-1$ if $\mu(x_i) < \infty$, and $a_i \in \mathbb{Z} \backslash \{0\}$ if $\mu(x_i)= \infty$. Since $g\neq 1$ we have $k\geq 1$. It follows from the previous lemma that $\rho_g(e_{x_1}) = n_{x_1}^{a_1} e_{x_1}$. Since $a_1 \in \mathbb{Z} \backslash \{0\}$ if $\mu(x_1) = \infty$ and $1 \leq a_1 \leq \mu(x_1) - 1$ if $\mu(x_1) < \infty$, we get that $n_{x_1}^{a_1} \neq 1$, hence $\rho_g(e_{x_1}) \neq e_{x_1}$, hence $\rho_g \neq \mathrm{id}_V$.
	\item Assume that $x_1^{a_1} x_2^{a_2} \cdots x_k^{a_k}$ and $y_1^{b_1} y_2^{b_2} \cdots y_\ell^{b_\ell}$ are two words for $g\in \mathrm{Tr}(\Gamma)$ with the required properties. It follows from the previous Lemma applied with the word $x_1^{a_1} x_2^{a_2} \cdots x_k^{a_k}$ that $$\{x_1, x_2, \dots, x_k\} = \{ y\in V(\Gamma) \ \vert \ \rho_g(e_y) = \lambda e_z~\text{for some~}z\in V(\Gamma)\text{~and~}\lambda \neq 1\}.$$ But the last set only depends on $g$, and applying the previous Lemma with the word $y_1^{b_1} y_2^{b_2} \cdots y_\ell^{b_\ell}$ thus yields $\{x_1, \dots, x_k\} = \{y_1, \dots, y_\ell\}$. Since both sets are totally ordered, we get $k = \ell$, and $x_i = y_i$ for all $i=1, \dots, k$. It remains to show that $a_i = b_i$ for all $i=1, \dots, k$. So let $i\in \{1, 2, \dots, k\}$. By the previous Lemma applied with the word $x_1^{a_1} x_2^{a_2} \cdots x_k^{a_k}$ we have that $\rho_g(e_{x_i}) = n_{x_i}^{a_i} e_z$ for some $z\in V(\Gamma)$. But applying it with the word $x_1^{b_1}\cdots x_k^{b_k}$ yields that $\rho_g(e_{x_i}) = n_{x_i}^{b_i} e_{z'}$ for some $z'\in V(\Gamma)$. We thus have $z=z'$ and $n_{x_i}^{a_i}= n_{x_i}^{b_i}$ hence, since $n_{x_i}$ is an indeterminate if $\mu(x_i) = \infty$ and $1 \leq a_i, b_i \leq \mu(x_i)-1$ otherwise, we get that $a_i = b_i$. 
	\item This follows immediately from the previous point: if there is $x\in V (\Gamma)$ such that $\mu(x) = \infty$, then $x$ has infinite order, while otherwise the set of normal forms is finite, and there are $\prod_{x\in V(\Gamma)} \mu(x)$ such forms. 
	\item If $\mu(x) = \infty$ for all $x\in V(\Gamma)$, by Remark~\ref{rmq_monoid} we can still write every element of $\mathrm{Tr}^+(\Gamma)$ under the form $x_1^{a_1} x_2^{a_2} \cdots x_k^{a_k}$, where $x_i \in V(\Gamma)$, $x_1 \succ x_2 \succ \cdots \succ x_k$, and $a_i \geq 1$ for all $i$. Since the image of such a word in $\mathrm{Tr}(\Gamma)$ is a normal form, distinct elements of $\mathrm{Tr}^+(\Gamma)$ cannot be identified inside $\mathrm{Tr}(\Gamma)$, because it would give two different normal forms for an element of $\mathrm{Tr}(\Gamma)$, contradicting the point $2$.
	\item Let $g\in \mathrm{Tr}(\Gamma)$ with $g\neq 1$. Write $g$ in the form $x_1^{a_1} x_2^{a_2} \cdots x_k^{a_k}$, where $x_i \in V(\Gamma)$, $x_1 \succ x_2 \succ \cdots \succ x_k$, and $a_i \in \mathbb{Z} \backslash \{0\}$. By the previous Lemma we have $\rho_g(e_{x_1}) = q_{[x_1]}^{a_1} e_{x_1}.$ For $n\geq 1$ we thus have $\rho_g^n (e_{x_1}) = q_{[x_1]}^{n a_1} e_{x_1}$. Since $a_1 \neq 0$ and $q_{[x_1]}$ is an indeterminate, we have $q_{[x_1]}^{n a_1} \neq 1$, hence $\rho_g^n (e_{x_1}) \neq e_{x_1}$, hence $\rho_g^n \neq \mathrm{id}_V$, hence $g^n \neq 1$.     
\end{enumerate}
\end{proof}

\begin{rmq}
In~\cite{trickle}, it is shown that when $\mu(x) = \infty$ for all $x\in V(\Gamma)$ and $V(\Gamma)$ is finite, the group $\mathrm{Tr}(\Gamma)$ is torsion-free. The authors ask whether this remains true when $V(\Gamma)$ is infinite. Point $(4)$ of the above theorem shows that this is the case when $\Gamma$ is complete.    
\end{rmq}

\begin{rmq}[Word problem]
If $\Gamma$ is a trickle graph with $V(\Gamma)$ is finite (and $\Gamma$ not necessarily complete), then by~\cite[Corollary 2.5]{trickle} the word problem is solvable in $\mathrm{Tr}(\Gamma)$. In the case where $\Gamma$ is complete, our results above give a new solution to the word problem: let $V(\Gamma)=\{x_1, x_2, \dots, x_n\}$, ordered so that if $x_i < x_j$, then $i > j$, and order the basis $\mathcal{B}=\{e_{x_1}, \dots, e_{x_n}\}$ accordingly. Let $g\in \mathrm{Tr}(\Gamma)$. Then Lemma~\ref{lem_main} together with Theorem~\ref{main}(2) tell us that the normal form of $g$ is given by $x_{i_1}^{a_1} x_{i_2}^{a_2} \dots x_{i_k}^{a_k}$, where $ 1 \leq i_1 < i_2 < \dots < i_k \leq n$ and the $x_{i_j}$ and their exponents $a_i$ are obtained as follows: the $i_j$'s are those integers for which the coefficient in the $i_j$'s column of the matrix $\mathrm{Mat}_{\mathcal{B}}(\rho_g)$ is different from $1$. For such an integer $i_j$, the coefficient is a nontrivial power of $n_{x_j}$, and this power (taken between $1$ and $\mu(x_{i_j})-1$  if $\mu(x_{i_j})<\infty$) is the exponent $a_j$. Since the matrix $\mathrm{Mat}_{\mathcal{B}}(\rho_g)$ can be calculated from any word for $g$, this gives an algorithm to determine the normal form of $g$.       
\end{rmq}

\section{An alternative proof of faithfulness for the trickle monoid}\label{sec_alt}

Again assume that $\Gamma$ is a complete trickle graph, and assume that $\mu(x) = \infty$ for all $x\in V(\Gamma)$. 

We give an alternative proof of the faithfulness of the representation $\rho$ restricted to $\mathrm{Tr}^+(\Gamma)$, using a philosophy that has been used to show that certain representations of Garside groups (such as braid groups~\cite{Krammer, Hee}, torus knot groups~\cite{Gobet_Burau}, etc.) are faithful: while Theorem~\ref{main} tells us that normal forms can be easily recovered from the representation, the strategy here is to show that those elements of $V(\Gamma)$ which left- or right-divide an element $g\in \mathrm{Tr}^+(\Gamma)$ can also be easily recovered from $\rho_g$.

We begin by recalling a few properties that we shall need. The following definition and lemma were given in~\cite{Hee} in the case of Artin monoids. 

\begin{definition}
	Let $M$ be a monoid admitting a homogeneous presentation $\langle S \mid R \rangle$. Let $N$ be a cancellative monoid and $\rho : M \longrightarrow N$ a morphism of monoids. For $m\in M$ let $I(m) = \{s\in S \mid s \leq_L m\}$, where $\leq_L$ denotes left-divisibility. We say that $\rho$ is \emph{initially injective} if for all $m,m'\in M$, $$ \rho(m) = \rho(m') \Rightarrow I(m) = I(m').$$  
\end{definition}

\begin{lemma}\label{lem_init_monoid}
	Assume that $M$ is a monoid admitting a homogeneous presentation $\langle S \mid R \rangle$. Let $N$ be a cancellative monoid and $\rho : M \longrightarrow N$ an initially injective morphism of monoids. Then $\rho$ is injective, and in particular, $M$ is cancellative.  
\end{lemma}

\begin{proof}
	Let $m, m'\in M$ such that $\rho(m) = \rho(m')$. Since the presentation $\langle S \mid R \rangle$ of $M$ is homogeneous, we have a well-defined length function $\ell =\ell_S$ on $M$. We argue by induction on $\ell(m) + \ell(m')$ to show that $m=m'$. If $\ell(m) + \ell(m') = 0$, then $m=m'=1$. Assume that $\ell(m) + \ell(m') > 0$. Then at least one element among $m$ and $m'$ is different from $1$, say $m \neq 1$. Let $s\in S$ such that $m = s x$ for some $x\in M$. Then $s\in I(m)$ and since $\rho$ is initially injective, we have $s \in I(m')$. Hence there exists $x'\in M$ such that $m' = s x'$. Now we have $$\rho(s) \rho(x) = \rho(sx) = \rho(m) = \rho(m') = \rho(s x') = \rho(s) \rho(x'),$$ and since $N$ is cancellative, we deduce that $\rho(x) = \rho(x')$. But $\ell(x) + \ell(x') = \ell(m) + \ell(m') - 2$, hence by induction we have $x = x'$, yielding $m = sx = sx' = m'$. We thus have that $\rho$ is injective, and since $N$ is cancellative, we deduce that $M$ is cancellative as well.   
\end{proof}

\begin{prop}\label{prop_des}
Let $m\in \mathrm{Tr}^+(\Gamma)$. Let $x\in V(\Gamma)$. Then $x \leq_L m$ if and only if there is $y\in V(\Gamma)$ such that $\rho_m(e_y) = \lambda e_x$ with $\lambda \neq 1$. 
\end{prop}

\begin{proof}
We argue by induction on $\ell(m)$. If $m=1$, then $\rho_m = \mathrm{id}_V$, and the claim is clear. Now assume that $\ell(m)>1$, and that the result holds by induction for elements $m'$ such that $\ell(m') < \ell(m)$. First suppose that $x \in V(\Gamma)$ is such that $x \leq_L m$. There is thus $m'\in \mathrm{Tr}^+(\Gamma)$ such that $m = xm'$. By definition of $\rho$, there is (a unique) $y\in V(\Gamma)$ such that $\rho_{m'}(e_y) = \lambda e_x$, for some $\lambda\in R$ of the form $\prod_{[x] \in O_\Gamma} q_{[x]}^{n_x}$, with $n_x \geq 0$. We then have $$\rho_g(e_y) = \rho_x \rho_{m'}( e_y) =\rho_x( \lambda e_x) = q_{[x]} \lambda e_x,$$ and $q_{[x]} \lambda \neq 1$ in $R$. 

Conversely, assume that there is $y\in V(\Gamma)$ such that $\rho_m(e_y) = \lambda e_x$ with $\lambda \neq 1$. Let $z\in V(\Gamma)$, $m'\in \mathrm{Tr}^+(\Gamma)$ such that $m = zm'$. If $z=x$ then $ x \leq_L m$ and we are done. We can thus assume that $z \neq x$. 

First assume that $z < x$. Then $\varphi_z(x)=x$, hence $\rho_z(e_x) = e_x$. We thus have that $$\rho_{m'} (e_y) = \rho_z^{-1} \rho_m(e_y) = \rho_z^{-1} \circ (\lambda e_x) = \lambda e_x,$$ hence by induction we have that $x \leq_L m'$. There is thus $m''\in \mathrm{Tr}^+(\Gamma)$ such that $m'= x m''$. Let $z' = \varphi_x^{-1} (z)$. We have $z' < x$ and thus the relation $\varphi_x(z') x = x z'$, that is, $zx = xz'$. We thus have $m = z m' = z x m'' = x z' m'',$ yielding $x \leq_L m$, as required. 

Now assume that $z \Vert x$. Then $\varphi_z(x) = x$, and $\rho_z( e_x) = e_x$. We thus have as in the previous case that $\rho_{m'}(e_y) = \lambda e_x$, and by induction, we have $x\leq_L m'$. There is thus $m''\in \mathrm{Tr}^+(\Gamma)$ such that $m'= x m''$. But in $\mathrm{Tr}^+(\Gamma)$ we have the relation $zx=xz$, hence $m = zm' = z x m'' = x z m'',$ again yielding that $x\leq_L m$. 

Finally, assume that $x < z$. Let $x'\in V(\Gamma)$ such that $\varphi_z(x') =x$. We have $\rho_{z}(e_{x'}) = e_x$. We thus have $$ \rho_{m'}(e_y) = \rho_{z}^{-1} \circ \rho_m( e_y) = \rho_z^{-1} (\lambda e_x) = \lambda e_{x'}.$$ 
By induction, we deduce that there is $m''\in \mathrm{Tr}^+(\Gamma)$ such that $m' = x' m''$. But $\varphi_{x'}(z) = z$ since $x' < z$, hence in $\mathrm{Tr}^+(\Gamma)$, we have the relation $zx'=xz$, and hence, $m = zm' = z x' m'' = xz m'',$ yielding $x\leq_L m$, as required. 
\end{proof}

\begin{cor}\label{cor_mon_inject}
The representation $\rho$ restricted to $\mathrm{Tr}^{+}(\Gamma)$ is injective, and the latter is cancellative. 
\end{cor}

\begin{proof}
Let $m, m'\in \mathrm{Tr}^+(\Gamma)$ such that $\rho_m = \rho_{m'}$. By Proposition~\ref{prop_des} we have $$I(m) = \{x \in V(\Gamma) \mid \exists y\in V(\Gamma), \exists \lambda\in R\backslash\{1\}\text{~such that~} \rho_m(e_y) = \lambda e_x\}= I(m').$$ It follows that $\rho$ is initially injective. We then apply Lemma~\ref{lem_init_monoid} to conclude the proof. 
\end{proof}

In the case where $V(\Gamma)$ is finite, by~\cite[Theorem 2.15]{trickle} the monoid $\mathrm{Tr}^+(\Gamma)$ is a Garside monoid. Another proof of the faithfulness of $\rho$ on the whole group $\mathrm{Tr}(\Gamma)$ can be deduced from this fact together with the corollary above, since in a Garside group, every element can be written as a fraction in two elements of the monoid. Given $g = m^{-1} m'$ with $m, m'\in \mathrm{Tr}^+(\Gamma)$ we thus have $$\rho(m^{-1} m') = \mathrm{id}_V \Leftrightarrow \rho(m') = \rho(m) \Leftrightarrow m = m' \Leftrightarrow g = 1.$$ 
Note that, if one does not want to use Theorem 2.15 from~\cite{trickle}, one can also argue as follows. It is enough in fact to show that $\mathrm{Tr}(M)$ is the group of fractions of $\mathrm{Tr}^+(M)$, and to this end, it suffices to show that $\mathrm{Tr}^+(M)$ is an Öre monoid, that is, a cancellative monoid in which every two elements admit a common left-multiple. Cancellativity was established in Corollary~\ref{cor_mon_inject} above. Letting $V(\Gamma) = \{x_1, x_2, \dots, x_n\}$ with $x_1 \succ x_2 \succ \dots \succ x_n$, let $x, y\in M$ and let $a_i, b_i \geq 0$ such that $x = x_1^{a_1} x_2^{a_2} \cdots x_n^{a_n}$ and $y = x_1^{b_1} x_2^{b_2} \cdots x_n^{b_n}$. It is not hard to check that $x_1^{c_1} x_2^{c_2} \cdots x_n^{c_n}$, where $c_i = \max(a_i, b_i)$, is a common left-multiple of $x$ and $y$. 
\begin{rmq}
It follows from Proposition~\ref{prop_des} together with Lemma~\ref{lem_main} that, in the notation of Lemma~\ref{lem_main}, for $m\in \mathrm{Tr}^{+}(M)$ written in the form $x_1^{a_1} x_2^{a_2} \cdots x_k^{a_k}$, we have $I(m) = \{ \varphi(x_i) \mid  i=1, \dots, k\}$. We recover in this specific case a result from~\cite{trickle}, namely Lemma~7.5. 
\end{rmq}

\begin{rmq}[Recovering atoms right-dividing an element]
For $m\in \mathrm{Tr}^+(\Gamma)$, define $R(m) = \{x\in V(\Gamma) \mid x \leq_R m\}$, where $\leq_R$ denotes right-divisibility. The opposite monoid $\mathrm{Tr}^+(\Gamma)^{\mathrm{op}}$ is also a trickle monoid with graph $\Gamma^{\mathrm{op}} = (\Gamma, \leq, \mu, (\varphi_x^{-1})_{x\in V(\Gamma)})$. We denote by $\rho^{\mathrm{op}}$ the representation from Section~\ref{sec_rep} associated with this opposite trickle group. Then for all $x\in V(\Gamma)$, we have that $\rho_x^{\mathrm{op}}$ is the transpose of $\rho_x$. It follows that for every $m\in \mathrm{Tr}^+(\Gamma)$, denoting $m^*$ its canonical image in $\mathrm{Tr}^+(\Gamma)^{\mathrm{op}}$, we have $\rho_{m^*}^{\mathrm{op}} = (\rho_m)^{t}$. For $m\in \mathrm{Tr}^+(\Gamma)$ we thus have 
\begin{align*}
R(m) &= \{x\in V(\Gamma) \mid x \leq_R m\} = \{x\in V(\Gamma) \mid x \leq_L m^*\} = I(m^*) \\ &=\{x \in V(\Gamma) \mid \exists y\in V(\Gamma), \lambda\in R\backslash\{1\}\text{~such that~} \rho_{m^*}^{\mathrm{op}}(e_y) = \lambda e_x\} \\ &=\{x \in V(\Gamma) \mid \exists y\in V(\Gamma), \lambda\in R\backslash\{1\}\text{~such that~} \rho_{m}^t(e_y) = \lambda e_x\} \\ &= \{x \in V(\Gamma) \mid \exists y\in V(\Gamma), \lambda\in R\backslash\{1\}\text{~such that~} \rho_{m}(e_x) = \lambda e_y\},
\end{align*}
where the fourth equality follows from Proposition~\ref{prop_des} applied with $\mathrm{Tr}^+(\Gamma)^{\mathrm{op}}$. We can thus also recover $R(m)$ from $\rho_m$. Moreover, writing $m$ in normal form $x_1^{a_1} x_2^{a_2} \cdots x_k^{a_k}$, it follows from Lemma~\ref{lem_main} that the last set above is equal to $\{x_1, x_2, \dots, x_k\}$, hence $R(m) =\{x_1, x_2, \dots, x_k\}$.  
\end{rmq}

	\end{document}